\documentclass[aop,MSNbibl,dvips]{arximspdf}

\doi{10.1214/12-AOP767} %kopijuoti is PTS
\volume{41}
\issue{3A}
\pubyear{2013}
\firstpage{1180}
\lastpage{1190}

\makeatletter
\def\rest{\upharpoonright}

\newtheorem{Theorem}{Theorem}[section]
\newtheorem{Proposition}[Theorem]{Proposition}
\newproclaim{Definition}[Theorem]{Definition}
\newproclaim{Remark}[Theorem]{Remark}
\newtheorem{Question}[Theorem]{Question}
\newtheorem{Lemma}[Theorem]{Lemma}
\def\E{\mathbf{E}}
\def\P{\mathbf{P}}
\def\Z{\mathbb{Z}}
\makeatother

\begin{document}
\begin{frontmatter}

\title{Invariant monotone coupling need not exist\thanksref{T1}}
\runtitle{Noninvariant Monotone Coupling}
\thankstext{T1}{Supported in part by
NSF Grants DMS-04-06017 and DMS-07-05518 and OTKA Grant K76099.}

\begin{aug}
\author[A]{\fnms{P\'eter} \snm{Mester}\corref{}\ead[label=e1]{pmester@indiana.edu}}%

\runauthor{P. Mester}
\affiliation{Indiana University}
\address[A]{Department of Mathematics\\
Indiana University\\
Bloomington, Indiana 47405-5701\\
USA\\
\printead{e1}}
\end{aug}

% HISTORY:
\received{\smonth{11} \syear{2010}}
\revised{\smonth{4} \syear{2012}}

% ABSTRACT
\begin{abstract}
We show by example that there is a
Cayley graph, having two invariant random subgraphs $X$ and $Y$,
such that there exists
a monotone
coupling between them in the sense that $X \subset Y$, although no such
coupling can be invariant. Here, ``invariant'' means that the
distribution is invariant under group multiplications.
\end{abstract}

% KEYWORDS


\end{frontmatter}

%s1 #&#
\section{Introduction}\label{s.bevezetes}
There are several models when one is selecting a
random subset of vertices or edges of a given graph $G=G(V,E)$
according to
some distribution. Formally these are $2^V$-valued random objects
where $V $ is the vertex set of $G$ (which can be replaced by $E$,
the set of edges). We can look at this as a $\{0,1\}$-labeling of the
vertices; then it is natural to allow more
general label sets ${\Lambda}$ replacing $\{0,1\}=2$.

We are interested in particular in Cayley graphs, and in this case,
most naturally occurring examples have an extra common feature:
{invariance}. This means that
their distribution is invariant under the group multiplication
of the base graph. More precisely, if $G$ is a right Cayley graph of
the group $\Gamma$, then the random object $\mathcal{ R}$ is
\textit{invariant} if for for any finite $\{v_1, \ldots, v_n\}\subset V$ and
$\gamma\in\Gamma$, the distribution of $(\mathcal{ R}(\gamma v_1),\ldots,
\mathcal{ R}(\gamma
v_n))$ does not depend on $\gamma$. Note that in this case $V=\Gamma$,
so it may seem confusing to use different notation. The reason is
that many concepts we define naturally generalize to the case of a
graph with a transitive group of automorphisms acting on the vertices,
and in general these are distinct notions.
The abundance of invariant processes on Cayley graphs motivates an
investigation of them in general. This was
done, for example, in {\cite{2}}.
% \ref b.BPLS1/.

In this context, our result is a counterexample. To explain it, we
first need to recall the notion of
{coupling}.

%de1.1 #&#
\begin{Definition}\label{d.coupling} If $\mathcal{ S}_1, \mathcal{
S}_2 $ are
random objects taking values in $ \Delta_1, \Delta_2 $, respectively,
then a \textit{coupling} of them is a random pair $(\tilde{ S}_1, { \tilde
S_2})$ taking values in $ \Delta_1 \times\Delta_2 $ such that for $i
\in\{1,2\}$, $\tilde{S}_i$ has the same distribution as $\mathcal{ S}_i$.
\end{Definition}

Intuitively this means that we
manage to produce the two objects using the same random source, so that
pointwise comparison makes sense.
Proofs using coupling arguments are usually very conceptual and fit
well with probabilistic intuition.

%re1.2 #&#
\begin{Remark}%\label{r.invariantcoupling}
If $\mathcal{ R}_i$ is a random $\Lambda_i$-labeling for $i\in\{1,2\}$,
then a coupling of them is a $\Lambda_1 \times\Lambda_2$-labeling,
and we say that this is an \textit{invariant coupling} if this labeling
is invariant in the above sense.
\end{Remark}

A very simple instance of this is the simultaneous coupling of
all\break
$\operatorname{Bernoulli}(E,p)$ percolations corresponding to possible
parameters $p \in[0,1]$
which we briefly recall.
A $\operatorname{Bernoulli}(E,p)$-percolation is obtained by putting i.i.d.
$\{0,1\}$-labels on the edges, where
for a given edge $e$, its label is $1$ with probability~$p$, and $0$
otherwise. Note that replacing the edge set $E$ with the set of
vertices $V$ in the above definition
is well defined, and we denote this process by $\operatorname{Bernoulli}(V,p)$.
There is a strong intuition that ``the bigger $p$ is, the bigger the
subgraph with label
$1$.''
We can make this intuition have a precise formal meaning as follows:
Put first i.i.d. uniform (from $[0,1]$) labels on the edges, which
we call $U$.
Then for each $p$, define a $\{0,1\}$-label $U_p$ so that if an edge has
$U$ label
$U(e)$, its $U_p$ label is $1$ if $U(e)\leq p$, and $0$ otherwise.
Clearly, as a distribution, $U_p$ is nothing but a $\operatorname{Bernoulli}(E,p)$-percolation, and
for $p\leq p^+$ we have $U_p \subset U_{p^+}$.

This is an example of what is called {monotone coupling}. For the
definition assume that the label set $\Lambda$ is partially ordered
by $\lesssim$.

%de1.3 #&#
\begin{Definition}\label{d.MonotoneCoupling} If $X$ and $Y$ are
random $\Lambda$-labelings of the same graph,
then we say that a coupling $(\tilde{ X}, \tilde{ Y})$ of $X$ and $Y$
is a \textit{monotone coupling} if $\tilde{ X}
\lesssim\tilde{ Y}$ almost surely.
\end{Definition}

The next two examples we mention are related to open questions which
motivates the question we are going to ask.

The first is the case of wired and free uniform
spanning forest measures ($\mathrm{WUSF}$ and $\mathrm{FUSF}$, resp.);
see {\cite{3}}.
These processes both can be considered as natural generalizations of
the uniform spanning tree (easily defined on finite graphs) to
infinite graphs.
It is known that there is a monotone coupling where the free one
dominates the wired one. However, in general, it is still open if there
is an invariant monotone coupling.

There are partial results which show that for certain classes of
graphs there indeed
exists an invariant monotone coupling between the $\mathrm{FUSF}$ and
$\mathrm{WUSF}$. For example, Lewis Bowen\vadjust{\goodbreak} {\cite{4}}
showed it for Cayley graphs of residually amenable groups, while
recently Russell Lyons and Andreas Thom {(personal communication, \cite{8})}
showed it
for the Cayley graphs of so-called \textit{sofic} groups.

The second example is random walk in random environment. In {\cite{1}},
Aldous and Lyons considered a continuous time
nearest-neighbor random walk $\operatorname{RW}(t,\mu)$, with jumps governed
by Poisson clocks on the edges with rates given by a distribution~$\mu$.
The walks start at the origin $o$ of the Cayley graph, and we are
interested in how different environments affect the return
probabilities.
In {\cite{1}}
they showed that if for two random environments $\mu_1,\mu_2$
(different clock frequencies in this case) there exists a monotone
coupling $\mu_1\leq\mu_2$, which is itself invariant, then
\[
\E_{\mu_1}\bigl(\P\bigl(\operatorname{RW}(t,\mu_1)=o\bigr)\bigr)
\geq\E_{\mu_2}\bigl(\P\bigl(\operatorname{RW}(t,\mu_2)=o\bigr)\bigr).
\]
We may ask what happens if we drop the condition for the coupling
being invariant. Is it enough, for example, that the marginals are
invariant?
Note that in \cite{1} they actually dealt with so-called
\textit{unimodular} processes, but this condition always holds for invariant
processes on Cayley graph (this fact is the mass transport principle
which we prove later).

Schramm and Lyons asked
(unpublished, {\cite{7}})
% \ref b.Oded-Russ/
the following; note that a positive answer would immediately settle
the above problems (note also that a more general question was asked in
\cite{1}, as Question 2.4):

%qu1.4 #&#
\begin{Question}\label{q.oded-russ} Let $X$ and $Y$ be invariant
subgraphs of a
Cayley graph~$\Gamma$, so that there exists a monotone coupling
between them. Does it follow that there exists a monotone coupling
between them which is also invariant?
\end{Question}

It is known that the answer to
the above question is ``yes'' if the Cayley graph is amenable; see
Proposition 8.6. in \cite{1}. In this paper
we show by an example that in full generality the answer is ``no.''

The Cayley graph we use is $T_3 \square C_n$ for $n$ large enough.
Here $T_3$ is the $3$-regular tree, and $C_n$ is the cycle of length
$n$, and, in general, for graphs $G$ and $H$ their \textit{Cartesian
product} $G \square H$ is the graph with vertex set $V(G\square
H)=V(G) \times V(H)$, and two vertices $(u_1,u_2)$ and $(v_1,v_2)$ are
connected in $G \square H$ if and only if either $u_1=v_1$ and $u_2$ is adjacent
with $v_2$ in $H$, or $u_2=v_2$ and $u_1$ is adjacent with $v_1$ in
$G$. It is easy to see that if $G_1,G_2$ are Cayley graphs of
$\Gamma_1,\Gamma_2$, respectively, then $G_1 \square G_2$ is a Cayley
graph of $\Gamma_1\times\Gamma_2$.

Note that $T_3$ is a Cayley graph of $\Z_2^{*3}:=\Z_2 * \Z_2 * \Z_2$
(here $H*K$ is the \textit{free product} of $H$ and $K$), and $C_n$ is a
Cayley-graph of $\Z_n$, so $T_3 \square C_n$ is a Cayley graph of
$\Z_2^{*3} \times\Z_n$.

For simplicity we make an assumption about $n$ which may not be optimal.
See Remark~\ref{r.nnagy} at the end of this section for an explanation.

%th1.5 #&#
\begin{Theorem}\label{t.alaptetel} If $n \geq376$, then there exist
two invariant
random $\{0,1\}$-labelings $X$ and $Y$ of $T_3 \square C_n$ so that
there is a coupling
$(\tilde{X},\tilde{ Y})$ of them for which $\tilde{ X} \leq
\tilde{ Y}$ holds, but no such coupling can be invariant.\vadjust{\goodbreak}
\end{Theorem}

The proof will be more succinct if we first show a similar result with
labels different from $\{0,1\}$. In this case the (partially ordered)
label set will be the power set $\mathcal{ P}(S)$ of some finite set $S$.
Note also that in this case we can use a tree as the underlying
Cayley graph:

%le1.6 #&#
\begin{Lemma}\label{l.alaplemma} If $n \geq376$ and $|S|=n$, then
there exist invariant
$\mathcal{ P}(S)$-labelings $\mathcal{ X}$ and $\mathcal{ Y}$ of $T_3$
so that there is a
coupling $(\tilde{\mathcal{ X}},\tilde{\mathcal{ Y}})$ of them for which
$\tilde{\mathcal{ X}} \subset\tilde{\mathcal{ Y}}$ holds, but no such
coupling can be invariant.
\end{Lemma}

Although the examples themselves might be artificial, they will have
some ``nice'' properties as well. So if we want to add some extra
conditions to Question~\ref{q.oded-russ} to get an affirmative answer, then
we know for sure that these nice properties will not work (at least
not alone). For a discussion of these, see the end of the last
section.

We summarize some conventions we use:
When $\mathcal{ S}$ is a random object
and $\mu$ is its distribution, we often will just express this by saying
that $\mathcal{ S}$ is a \textit{copy} of $\mu$, and in a similar way with a
further abuse of notation, if
$\mathcal{ T}$ is a random object with the same distribution as
$\mathcal{
S}$, we will also say that $\mathcal{ S}$ is a copy of $\mathcal{ T}$. We
also note that one way to specify a probability measure is to
describe a random object with the given measure as distribution. We
will do it
without further comments.

If the graph $G$ is understood, $V(G) $ will be its set of vertices
and $E(G)$ its set of edges. {We use right Cayley graphs and then left
multiplications are graph automorphisms.}

%re1.7 #&#
\begin{Remark}\label{r.nnagy} The condition that $n\geq376$ we made in
Theorem~\ref{t.alaptetel} was meant to ensure the following:
If $S$ is a finite set of cardinality $n$, and
$\alpha_1,\alpha_2,\ldots,  \alpha_{20}$ are i.i.d. uniform elements of
$S$, and $\beta_1,\ldots, \beta_9$ are also i.i.d. uniform elements of
$S$ (we emphasize that we make no extra assumption on the joint
distribution of the full family $\alpha_1,\ldots,
\alpha_{20},\beta_1,\ldots,\beta_9$), then with probability strictly
greater than ${1
\over2}$
the random elements $\alpha_1,\ldots, \alpha_{20}$ are all distinct, and
the random elements $\beta_1,\ldots,\beta_9$ are all distinct as well
(but it may happen that some $\beta_i=\alpha_j$). It is easy to see
that if $(1- \prod_{i=1}^{19}(1-{i\over n}))+(1-\prod_{i=1}^{8}(1-{i\over n}))< {1\over2}$ (which is
true for $n\geq376$), then this holds.
\end{Remark}

%s2 #&#
\section{The mass-transport principle and ends}\label{s.ends}

This section owes a lot to the exposition in
{\cite{6}}.
% \ref b.LP:book/.
An effective tool in showing that there is no invariant random process
on a Cayley graph satisfying a certain requirement is the so-called
mass-transport principle. {Recall that ${\Lambda}$ is the label set,
which will always be finite in our case, and $\Gamma$ is the group to
which the Cayley graph is associated.}
{The ``space of configurations'' $\Omega:={\Lambda}^{V}$
% or ${\Lambda}^{E}$
will be naturally equipped with the product ${\sigma}$-algebra.
Assume that $\mathcal{ R}$ is a probability measure on $\Omega$.
% so that
% ${\Omega}$ is either ${\Lambda}^V$ or ${\Lambda}^E$.
Note that $\Gamma$ acts on $\Omega$: for $\omega\in\Omega$, $\gamma
\in\Gamma$ and $v \in V$, let ${\gamma}{\omega}$ be the element of
$\Omega$ for which
${\gamma}{\omega}(v)={\omega}({\gamma^{-1}(v)})$.\vadjust{\goodbreak}

Let $F\dvtx  V\times
V\times\Omega\rightarrow[0, \infty]$ be a diagonally invariant
measurable function [meaning that $F(x,y, \omega)=F(\gamma x, \gamma
y, \gamma\omega)$ for all $\gamma\in\Gamma$].
The quantity $F(x,y, \omega)$ is often called \textit{the mass sent by
$x$ to $y$} or \textit{the mass received by $y$ from $x$}, and then $F$ is
thought to describe a ``mass transport'' among the vertices which may
depend on some randomness created by $\mathcal{ R}$.
The mass-transport principle says that if $\mathcal{ R}$ is invariant,
then for the identity $o \in V$ the expected overall mass $o$ receives
is the same as the expected overall mass it sends out. Now we
formalize and prove this:

%th2.1 #&#
\begin{Theorem}\label{t.masstrans}
If $\mathcal{ R}$ and $F$ are as above, $\mathcal{ R}$ is invariant,
$f(x,y):=\E_\mathcal{ R} F(x,y,*)$, then

\[
\sum_{x\in V}f(o,x)=\sum_{x\in V}f(x,o).
\]
\end{Theorem}

To prove it, first observe that the invariance of $\mathcal{ R}$ implies
that $f$ is also diagonally invariant. This implies that
$f(o,x)=f(x^{-1}o,x^{-1}x)=f(x^{-1},o)$, and this finishes the proof
since inversion is a bijection.

% Note that the same result and
%proof goes through if we consider
%the space of configurations to be $\Lambda^E$ instead.
%However, we will not need that form
%(although we will use edge labeling as an auxiliary tool).

This means that in order to show that a random process with a given
property cannot be invariant, it is enough to show that the property
in question allows us to define a mass transport contradicting the
above equality. We emphasize that it is important here that we mean
invariant processes on a Cayley graph and not just on a graph which
has a transitive group of automorphism. The notion of an ``end'' in a
tree, which we are about to define, will also lead to an example where
the obvious generalization of the mass transport principle to an
arbitrary transitive graph fails.

% \textbf{as we will do at the end of \ref s.lalley/.}

From now until Section~\ref{s.szam}, the base graph is always $T_3$, the
$3$-regular tree.
If $v$ is a vertex, then $J(v)$ will
denote the set of edges for which $v$ is one of the endpoints.

A
``ray'' is a one-sided
infinite path (i.e., a sequence of vertices $v_0, \ldots, v_n, \ldots$
so that there is no repetition and $v_i$ and $v_{i+1}$ are
adjacent).
We call two rays equivalent if their symmetric difference is
finite. An equivalence class is then called an \textit{end}.
If we fix an end $\xi$, then for any vertex $v$ there is
a unique ray $ v =v_0^{\xi}, v_1^{\xi}, \ldots, v_n^{\xi}, \ldots $
so that the ray
starts at $v$ and belongs to the
equivalence class $\xi$. Let the unique
edge joining $v$ with $v_1^{\xi}$ be $e_{v \rightarrow\xi}$, and
let us denote
$J(v)-\{e_{v \rightarrow\xi}\}$ as $J^{\xi}(v)$. Observe that for
distinct vertices $v_1,v_2$, we have
\[
J^{\xi}(v_1)\cap J^{\xi}(v_2)=
\varnothing.
\]

This will be important in constructing a monotone coupling of our
processes, and it also implies that an end cannot be determined using
invariant processes. The intuition is simple: given an end
$\xi(\omega)$ (which is ``somehow determined'' by a configuration
$\omega$) a vertex $v$ could send mass $1$ to each of the two
vertices that are the other endpoints of the two edges in
$J^{\xi(\omega)}(v)$. In this way the overall mass\vadjust{\goodbreak} sent out is $2$,
while the overall mass received is $1$.
To make this precise in a general setting, we have to deal with
measurability issues related to how a configuration~${\omega}$
determines an end $\xi({\omega})$, but this is not needed for our
purposes. While it is not important for our later work, note that if
we put extra edges into $T_3$ by connecting every vertex $v$ with~$v_2^{\xi}$, then we get a transitive graph where the obvious
generalization of the mass transport principle fails.

%s3 #&#
\section{The fixed-end trick}\label{s.lalley}

As we have indicated, an end cannot be determined using invariant
processes in a tree, and Lalley (unpublished, {\cite{5}})
proposed a way to exploit this fact to settle Question~\ref{q.oded-russ}.
Here we present a simpler version of the idea; see the last paragraph
in this section for the original one.
Given an end $\xi$ in $T_3$, we shall
define a $\{0,1\}^2$-labeling
$(\mathrm{X}^{\xi}, \mathrm{Y}^{\xi})$ so
that its components,
$\mathrm{X}^{\xi}$ and $\mathrm{Y}^{\xi}$, are invariant,
and $(\mathrm{X}^{\xi}, \mathrm{Y}^{\xi})$ is a monotone
coupling of them, that is,
\[
\mathrm{X}^{\xi} \leq \mathrm{Y}^{\xi}.
\]

Let $\{ \eta
(e) \}_{e \in E}$ be a $\operatorname{Bernoulli}(E,{1 \over2})$ label.
For a vertex $v$, let
\[
\mathrm{X}^{\xi}(v):= \operatorname{ max} \bigl\{ \eta(e); e \in
J^{\xi}(v) \bigr\},
\]
while
\[
\mathrm{Y}^{\xi}(v):=\operatorname{ max } \bigl\{\eta(e); e\in J(v) \bigr\}.
\]

It is clear that $\mathrm{Y}^{\xi}$ itself is an invariant labeling.

However, $\mathrm{X}^{\xi}$
is also invariant since the family $\mathrm{X}^{\xi}(v)_{v \in V}$
is actually i.i.d.! This is because of the observation from the last
section that $J^ {\xi} (v_1)$ and $J^ {\xi} (v_2)$ are disjoint for
$v_1 \not = v_2$.
So
$\mathrm{X}^{\xi}$ itself is actually $\operatorname{Bernoulli}(V, {3 \over4})$.
Since the monotone coupling of these processes was defined using an
end (a noninvariant step), it is reasonable that maybe these
processes already witness Theorem~\ref{t.alaptetel}.

However, the construction below---which is due to Peres
(unpublished,~{\cite{9}})---shows that there exists an invariant monotone coupling between
$\mathrm{X}^{\xi}$ and $\mathrm{Y}^{\xi}$.

%pr3.1 #&#
\begin{Proposition}\label{p.yuval}
Let $\{ \eta
(e) \}_{e \in E}$ be as above.
For each vertex $v$ with $J(v)=:\{e_1(v),e_2(v),e_3(v)\}$,
define $\hat{\mathrm{X}}(v):=0$ if and only if $\{ \eta(e_1)= \eta(e_2)
=\eta(e_3)\}$, and $\hat{\mathrm{X}}(v):= 1$ otherwise.

Then $(\hat{\mathrm{X}},\mathrm{Y}^{\xi})$ is an invariant and monotone
coupling of $\mathrm{X}^{\xi}$ and $\mathrm{Y}^{\xi}$.
\end{Proposition}

\begin{pf}
It is clear that $\hat{\mathrm{X}}\leq\mathrm{Y}^{\xi}$ and the coupling
$(\hat{
\mathrm{X}}, \mathrm{Y}^{\xi})$ is clearly invariant.
What we need to show is that the above defined $\hat{\mathrm{X}}$ is a
$\operatorname{Bernoulli}(V,{3\over
4})$ vertex labeling (i.e., a copy of $\mathrm{X}^{\xi}$) so $(\hat{\mathrm{X}}, \mathrm{Y}^{\xi})$ is a coupling of $\mathrm{X}^{\xi}$ and $\mathrm{Y}^{\xi}$.

First let us introduce a notation: if $V_1,V_2$ are finite disjoint
sets of vertices, then let $\mathrm{S}[V_1,V_2]:= \{ \hat{\mathrm{X}} \rest
V_1 =1, \hat{\mathrm{X}} \rest V_2 =0\}$.
We show that $\hat{\mathrm{X}}$ is $\operatorname{Bernoulli}(V,{3\over4})$
directly by proving that $\P( \mathrm{S}[V_1,V_2])=({3 \over
4})^{|V_1|}({1 \over4})^{|V_2|}$.\vadjust{\goodbreak}

The proof goes by induction on $|V_1 \cup V_2|$. The statement holds
when $|V_1 \cup V_2|=1$.

To proceed, consider the subgraph spanned by $V_1 \cup V_2$. This is a
forest, so it has some vertex ${t}$ which is either a leaf or an
isolated point (i.e., $t$ has at most one neighbor in $V_1 \cup V_2$).
Let $e_1,e_2, e_3$ be the edges emanating from~$t$, and assume that the
other endpoints of $e_1,e_2$ are \textit{not} in $V_1 \cup V_2$.

A key observation is that ``flipping'' the $\eta$ labels of each edge
leaves the $\hat{\mathrm{X}}$ labels unchanged. That means that for any
pair of finite disjoint vertex sets $W_1,W_2$ and \textit{any} edge $e$
the event
$\{\eta(e)=1\}$ [and similarly $\{\eta(e)=0\}$] cuts $ \mathrm{S}[ W_1,
W_2]$ exactly in half:
$\P(\{\eta(e)=1\}\cap \mathrm{S}[ (W_1, W_2])=P(\{\eta(e)=0\}\cap \mathrm{S}[ W_1, W_2]) = {1 \over2}P( \mathrm{S}[ W_1, W_2])$.

Using the $\eta$-flipping observation it is enough to show that
$\P(\{\eta(e_3)=1\} \cap\mathrm{S}[V_1,V_2])=({1 \over2})({3 \over
4})^{|V_1|}({1 \over4})^{|V_2|}$.

Consider first the case where $t\in V_1$.
In that case $\{\eta(e_3)=1\} \cap\mathrm{S}[V_1,V_2]=\{\hat{\mathrm{X}}(t)=1\} \cap\{\eta(e_3)=1\}
\cap\mathrm{S}[V_1-\{t\},V_2]=\{\mbox{at least one of }\eta(e_1)\mbox{ and }\eta(e_2)\break\mbox{ is } 0\} \cap
\{\eta(e_3)=1\} \cap\mathrm{S}[V_1-\{t\},V_2]$. The point of this is
that $\{\mbox{at least one of}\break\eta(e_1)\mbox{ and }\eta(e_2)\mbox{ is }0\}$ and
$\{\eta(e_3)=1\} \cap\mathrm{S}[V_1-\{t\},V_2]$ are independent,
and by induction and the $\eta$-flipping observation, we know the
probability of the latter; it is $({1 \over2})({3 \over
4})^{|V_1|-1}({1 \over4})^{|V_2|}$. Combining this with the fact that
$\P\{\mbox{at least}\break \mbox{one of }\eta(e_1)\mbox{ and } \eta(e_2)\mbox{ is }0 \}={3
\over4}$ gives us $\P(\{\eta(e_3)=1\} \cap\mathrm{S}[V_1,V_2])=\break({1
\over2})({3 \over4})^{|V_1|}({1 \over4})^{|V_2|}$.

If $t\in V_2$ we have $\{\eta(e_3)=1\} \cap\mathrm{S}[V_1,V_2]=
\{\hat{\mathrm{X}}(t)=0\} \cap\{\eta(e_3)=1\}\cap\mathrm{S}[V_1,V_2-\{t\}]=\{ \eta(e_1)=\eta(e_2)=1 \} \cap\{\eta(e_3)=1\}
\cap\mathrm{S}[V_1,V_2-\{t\}]$. The independence of $ \{
\eta(e_1)=\eta(e_2)=1 \}$ and $\{\eta(e_3)=1\} \cap\mathrm{S}[V_1,V_2-\{t\}]$ combined with what we know by induction proves the
claim again.\vspace*{-2pt}
\end{pf}

Although the above processes could be coupled in an invariant way, it
is clear that the idea leaves us a lot of freedom to use other
partially ordered sets and other monotone operations (instead of
taking maxima, we could take the sum, e.g., which was Lalley's
original suggestion). But it seems that other examples are difficult
to analyze from the point of view of Question~\ref{q.oded-russ}. With our next
construction, however, it will be very succinct why a monotone
coupling cannot be invariant.\vspace*{-2pt}

%s4 #&#
\section{Set valued labels on $T_3$}\label{s.halmazok}

In this section, we describe an example that will prove Lemma~\ref{l.alaplemma}.

Let $S$ be a finite set with $|S|=n\geq376$, and let $\mathcal{ P}(S)$
denote its power set. We will use $\mathcal{ P}(S)$ as a label set with
inclusion as a partial order. The two invariant $\mathcal{
P}(S)$-labelings $\mathcal{ Y}_S$ and $\mathcal{ X}_S$ of the vertices of
$T_3$ are defined as follows (for the rest of this section we drop the
subscript $S$ but in the next section we use it again).

To construct $\mathcal{ Y}$, we first label the edges of $T_3$ with
independent uniform elements from ${S}$. Let us call this labeling
$\lambda$. Then for a vertex $v$, let $\mathcal{ Y}(v):=\bigcup_{e \in
J(v)} \{\lambda(e)\}.$\vadjust{\goodbreak}

To construct $\mathcal{ X}$, we first define its marginal $\nu$ on the
vertices. To get a copy of $\nu$ first pick a uniform $(x_1,x_2)\in
S\times S$, and then take $\{x_1\} \cup\{x_2\}$. Finally let $\{\mathcal{
X}(v)\}_{v \in V(T_3)}$ be a labeling of the vertices with i.i.d.
copies of ${\nu}$.

%re4.1 #&#
\begin{Remark}\label{r.ymetszet} Observe that if $\hat{\mathcal{ Y}}$
is any copy of
$\mathcal{ Y}$, then the following is true: if $v_0$ is any vertex with
neighbors $v_1,v_2,v_3$, then any $s \in\hat{\mathcal{ Y}}(v_0)$ is also
contained in at least one of the $\hat{\mathcal{ Y}}(v_i)$'s for $i\in
\{1,2,3\}$.
\end{Remark}

By fixing an end $\xi$, we can present a monotone coupling of $\mathcal{
X}$ and $\mathcal{ Y}$ just as in Lalley's example. Get the copy of
$\mathcal{ Y}$ in the exact same way as above using $\lambda$ as a source,
but also use this $\lambda$ to get the copy of $\mathcal{ X}$ as
$ \mathcal{ X}^{\xi}(v):= \bigcup_ {e \in J^{\xi}(v)} \{ \lambda(e) \}$.
Then clearly
$\mathcal{ X}^{\xi}(v)\subset\mathcal{ Y}(v)$ holds for all $v$, and
$\mathcal{
X}^{\xi}$ is indeed a copy of $\mathcal{ X}$ [recall that the disjointness
of the different $J^{\xi}(v)$'s guarantees the independence for
different vertices and the marginals are clearly the same].

However there cannot be any invariant monotone coupling as we
will show now (which together with the previous paragraph proves Lemma~\ref
{l.alaplemma}).

%pr4.2 #&#
\begin{Proposition}\label{p.alaplemma1}
There exists no coupling of $\mathcal{ X}$ and $\mathcal{ Y}$ which is
both invariant and monotone.

\end{Proposition}

\begin{pf}
Let $(\mathcal{ X}^*,\mathcal{ Y}^*)$ be any monotone coupling
of $\mathcal{ X}$ and $\mathcal{ Y}$. We will show that using this monotone
coupling, we can define a mass transport $F$ which contradicts the mass
transport principle, showing that the coupling cannot be invariant.
To define the mass transport we have to say for every pair $(v_0,v)$
of vertices and every possible configuration $\omega$ [defined in
terms of $(\mathcal{ X}^*,\mathcal{ Y}^*)$] the value $F(v_0,v,\omega
)\in
[0,\infty]$. The dependence on $\omega$ will be through an event
$E(v_0)$ which we define now.

%de4.3 #&#
\begin{Definition}\label{d.nnagy}
First, let $v_1,v_2,v_3$ be the neighbors of $v_0$ and $v_4,\ldots,v_9$
be the vertices at graph distance $2$ from $v_0$ (in any order).

We say that $E_1(v_0)$ holds if for each $1\leq i, j\leq3, i\not
=j$, we have $|\mathcal{ Y}^*(v_j)|=3$, $\mathcal{ Y}^*(v_j) \cap
\mathcal{
Y}^*(v_i)=\varnothing$.

We say that $E_2(v_0)$ holds if for each
$0 \leq i,j \leq 9, i\not=j$, we have $|\mathcal{ X}^*(v_i)|=2$ and
$\mathcal{ X}^*(v_i)\cap\mathcal{ X}^*(v_j)=\varnothing$.

Finally, let $E(v_0):=E_1(v_0)\cap E_2(v_0)$.
\end{Definition}

Note the connection with the condition in Remark~\ref{r.nnagy}: the labels
$\mathcal{ X}(v_0), \ldots, \break \mathcal{ X}(v_9)$ can be identified with
$\{\alpha_1,\alpha_2\}, \ldots, \{\alpha_{19},\alpha_{20}\},$ while the
edge labels of those $9$ edges which are relevant in the $\mathcal{ Y}$
labels of $v_0,v_1,v_2,v_3$ can be identified
with $\beta_1,\ldots, \beta_9$. Then the condition we made on $n$
ensures that $\P(E(v_0))> {1\over2}$.

Now we are ready to define the mass transport $F\dvtx  V \times V \times
\Omega\rightarrow[0,\infty]$. If $E(v_0)$ does not hold, then set
$F(v_0,v, \omega):=0$ for each vertex $v$.\vadjust{\goodbreak}
If $E(v_0)$ holds, then let $F(v_0, v, \omega):=1$ if $v_0$ is a
neighbor of $v$ and $\mathcal{ X}^*(v_0) \cap\mathcal{ Y}^*(v) \not=
\varnothing$, while in every other case, set $F(v_0,v,\omega):=0$.

We claim that the expected mass the origin sends out is strictly
greater than
${1}$, while the mass it receives is not greater than $1$ (even
point-wise).

To prove this we show first that if $E(v_0)$ holds, then the mass
$v_0$ sends out is \textit{exactly} $2$. This combined with the fact that
$\P(E(v_0))> {1 \over2}$ implies the first part of the claim.
Let $\mathcal{ X}^*(v_0)=:\{s_1,s_2\}$; observe that $E_2(v_0)$ implies
$s_1 \not=s_2$. By monotonicity of the coupling, $\{s_1,s_2\}\subset
\mathcal{ Y}^*(v_0)$, so by Remark~\ref{r.ymetszet} there exist neighbors
$v_0(s_1),v_0(s_2)$ of $v_0$ so that $\mathcal{ Y}^*(v_0(s_i))$ contains
$s_i$. Since $\mathcal{ Y}^*(v_1), \mathcal{ Y}^*(v_2),\mathcal{
Y}^*(v_3)$ are
pairwise disjoint sets [by $E_1(v_0)$], there can be at most two of
them which nontrivially intersect $\{s_1,s_2\}$, and this implies
$v_0(s_1)\not=v_0(s_2)$. By definition, $v_0(s_1)$ and $v_0(s_2)$ are
exactly the vertices receiving nonzero mass from $v_0$.

To prove that the expected mass $v_0$ receives is at most $1$, assume
that $v_0$ receives nonzero mass from $v_1$ and $v_2$. First of all,
$v_1$ sends out nonzero mass only if $E(v_1)$ [and in particular
$E_2(v_1)$] holds.
Since $v_0$ and $v_2$ are
both within distance~$2$ from $v_1$, the event $E_2(v_1)$ implies that
$\{a_1,a_2\}:=\mathcal{ X}^*(v_1),\{b_1,b_2\}:=\mathcal{ X}^*(v_2)$, and
$\{c_1,c_2\}:=\mathcal{ X}^*(v_0)$ are pairwise disjoint and each has size
$2$. By the condition for the mass transport, $\mathcal{ Y}^*(v_0)$
contains one of the $a_i$'s, one of the $b_i$'s and---by the
monotonicity of the coupling---$\{c_1,c_2\}$ as well. But this would
mean that $\mathcal{ Y}^*(v_0)$ has at least four distinct elements,
which is impossible.

This mass transport violates the mass transport theorem,
so no monotone coupling of $\mathcal{ X}$ and $\mathcal{ Y}$ can be invariant.
\end{pf}

%re4.4 #&#
\begin{Remark}\label{r.endlike} Observe that an end $\xi$ can be
identified by the
orientation on the edges given as follows: orient the edges in
$J^{\xi}(v)$ \textit{away} from $v$. Then the out-degree of a vertex is
always $2$ while the in-degree is always $1$. The mass transport above
has some similarity with this end: if for a vertex $v$ we define
$J^{(\mathcal{ X}^*,\mathcal{ Y}^*)}(v)$ to be the set of edges
connecting $v$
with vertices receiving nonzero mass from $v$ and we orient the edges
in $J^{(\mathcal{ X}^*,\mathcal{ Y}^*)}(v)$ away from $v$, then the out-degree
of a vertex is either $0$ or $2$ and the in-degree is either $0$ or
$1$.\looseness=-1
\end{Remark}

%s5 #&#
\section{\texorpdfstring{The $\{0,1\}$-labels on $T_3 \Box C_n$}{The $\{0,1\}$-labels on $T_3 Box C_n$}}\label{s.szam}

Now we prove Theorem~\ref{t.alaptetel}. The Cayley graph we use is $T_3
\square C_n$ defined in the \hyperref[s.bevezetes]{Introduction.} Recall that $T_3 \square
C_n$ is a Cayley graph of $\Z_2^{*3} \times\Z_n$. If we want to check
the invariance of a process defined on $T_3\square C_n$, it is enough
to check invariance under group multiplication from
$\Z_2^{*3}$ and $\Z_n$ since the direct components generate the full group.

The processes we define can be considered as very faithful copying of
the previous processes. In the previous section, the set $S$ whose
subsets were used as labels was not important besides its cardinality
$|S|$. Now it will be convenient to choose it to be $S:=V(C_n)$.
If $\mathcal{ Z}$ is any $\mathcal{ P}(S)$-labeling of the vertices of $T_3$,
then let $\operatorname{lift}(\mathcal{ Z})$ be the following $\{0,1\}$-labeling of
$T_3 \square C_n$:
for a vertex $(u,v) \in V(T_3 \square C_n)$, let $\operatorname{lift}(\mathcal{
Z}) (u,v):=1$ if $v \in\mathcal{ Z}(u)$, otherwise let $\operatorname{lift}(\mathcal{
Z})(u,v):=0$. Note that this function $\operatorname{lift}$ from the $\mathcal{
P}(S)$-labelings of $T_3$ to the $\{0,1\}$-labelings of $T_3 \square
C_n$ is invertible.

Consider the previously defined processes $\mathcal{ X}_S,\mathcal{ Y}_S$.
Let $X:=\operatorname{lift}(\mathcal{ X}_S)$ and $Y:=\operatorname{lift}(\mathcal{ Y}_S)$.

%pr5.1 #&#
\begin{Proposition}\label{p.alaptetel1}
The above defined $X$ and $Y$ witness the truth of Theorem~\ref{t.alaptetel}.
\end{Proposition}

\begin{pf}
First, the invariance of $X$ and $Y$ under group multiplication from
$\Z^{*3}_2$ follows from the fact that $\mathcal{ X}_S$ and $\mathcal{ Y}_S$
were invariant on $T_3$, and the invariance under $\Z_n$ follows from
the fact that for a fixed vertex $v_0$, the distribution of $\mathcal{
X}_S$ ($\mathcal{ Y}_S$) is invariant under any permutation of $S$.

Second, there exists a monotone coupling of $X,Y$ since if $(\mathcal{
X}^*_S,\mathcal{ Y}^*_S)$ is any monotone coupling of $\mathcal{
X}_S,\mathcal{
Y}_S$, then $(\operatorname{lift}(\mathcal{ X}^*_S),\operatorname{lift}(\mathcal{
Y}^*_S))$ is
clearly a monotone coupling of $X$ and $Y$.

Third, if $(X^*,Y^*)$ was an invariant monotone coupling, then
$(\operatorname{lift}^{-1}(\mathcal{ X}^*),\break  \operatorname{lift}^{-1}(\mathcal{ Y}^*))$ would
have been an invariant coupling of $\mathcal{ X}_S,\mathcal{ Y}_S$,
which is
impossible as we have seen.
\end{pf}

It would be nice to have some natural condition on random subgraphs
under which the answer to Question~\ref{q.oded-russ} would be affirmative.
Here we point out two conditions which are ruled out by our example.

A random subgraph $Z$ is said to be $k$-\textit{dependent} if for vertex
sets $S_1,S_2,\ldots,\break S_m$ whose pairwise distances are all at least
$k$, the random objects $F_i:=Z \rest S_i, 1\leq i \leq m,$ are
independent.
Our example is $k$-dependent for large enough~$k$ (depending on the
cycle size $n$).
So assuming $k$-dependence is certainly not enough.

With slight modifications, we can exclude other conditions as well.
Observe that the mass transport we used would still work [in the sense
that $E(v)$ would have probability greater than ${1\over2}$] if we
``perturbed'' our processes with a $\operatorname{Bernoulli}(V, {\varepsilon})$
process for $\varepsilon>0$ small enough [meaning that we change the
original labels on those vertices where $\operatorname{Bernoulli}(V,
{\varepsilon})$ turns out to be $1$].
A random subgraph $Z$ is said to have \textit{uniform finite energy} if
there exists an $\varepsilon\in(0,1)$ so that for a vertex~$v$ we have
$\varepsilon<\P(Z(v)=1|Z \rest V-\{v\})<1-\varepsilon$.
By using this idea of perturbing the labels we see that assuming that
the process has uniform finite energy is not enough either.

\section*{Acknowledgments}
The author thanks Russell Lyons for suggesting working on this problem
and explaining some background, in particular Lalley's idea.
The timeline was reconstructed with the kind help of Russell Lyons and
Yuval Peres: after the question was asked in 1997,
the problem was discussed in a meeting in 1998 where the ideas included
in our ``unpublished'' references have arisen.

I thank the referee for careful reading and for useful comments.

% imsref loaded by akundreckaite, 2012-10-17 12:59:41

\printaddresses

\end{document}